\documentclass[a4paper,
               aps,
               prl,
               amsfonts,
               amssymb,
               amsmath,
               nobibnotes,
               twocolumn, 
               twoside,
               balancelastpage,
               eqsecnum] 
               {revtex4-1}
\usepackage[utf8]{inputenc} 
\usepackage[english]{babel} 

\usepackage{amsmath, amssymb}
\usepackage{url}
\usepackage{wrapfig}
\usepackage{subfiles}
\usepackage[squaren]{SIunits}
\usepackage{caption}
\usepackage[normalem]{ulem}
\usepackage{subfigure}
\usepackage{lipsum}
\usepackage{todonotes} 

\begin{document}
\selectlanguage{english}
\bibliographystyle{apsrev4-1}

\setcounter{section}{1} 
\setcounter{equation}{0} 

\author{Arianna Marchionne}
\author{Peter Ditlevsen}
\affiliation{Niels Bohr Institute, University of Copenhagen, Denmark}
\author{Sebastian Wieczorek}
\affiliation{Department of Applied Mathematics, University College
  Cork, Ireland}

\date{\today}

\title{Three types of nonlinear resonances 
}

\begin{abstract}
  We analyse different types of nonlinear resonances in a weakly
  damped Duffing oscillator using bifurcation theory techniques. In
  addition to
  (i) odd subharmonic resonances found on the primary branch of
  symmetric periodic solutions with the forcing frequency and (ii)
  even subharmonic resonances due to symmetry-broken periodic
  solutions that bifurcate off the primary branch and also oscillate
  at the forcing frequency, we uncover (iii) novel resonance type due
  to isolas of periodic solutions that are not connected to the
  primary branch. These occur between odd and even resonances,
  oscillate at a fraction of the forcing frequency, and give rise to a
  complicated resonance `curve' with disconnected elements and high
  degree of multistability.

  We use bifurcation continuation to compute resonance tongues in the
  plane of the forcing frequency vs. the forcing amplitude for
  different but fixed values of the damping rate. In this way, we
  demonstrate that identified here isolated resonances explain the
  intriguing structure of ``patchy tongues'' observed for week
  damping and link it to a seemingly unrelated phenomenon of
  ``bifurcation superstructure'' described for moderate damping.

\end{abstract}

 \maketitle

\section{Introduction}
The response of non-linear oscillators to external periodic forcing is
a surprisingly complex phenomenon. It has been investigated in a whole
range of natural systems ranging from
biology~\cite{FitzHugh:1961,Nagumo:1962} to glacial cycles
~\cite{Crucifix:2012, Saedeleer:2013}. Also, various physical systems,
such as lasers~\cite{Wieczorek:2005}, electric
circuits~\cite{van_der_Pol:1926} and mechanical
oscillators~\cite{Duffing:1918} have been studied through the last
century.  A widely studied model system of a \emph{damped and forced
  nonlinear oscillator} is the Duffing oscillator, which is an
anharmonic oscillator that can be realised as a mass on a spring with
a non-linear restoring force~\cite{Duffing:1918, Holmes:1976, Holmes:1979, Guckenheimer:2013}:
\begin{equation}
\frac{d^2x}{dt^2} + \gamma\frac{dx}{dt}  + \omega_0^2 x + \beta x^3 = A \cos(\omega t).
\label{eqn:duffing}
\end{equation}
Here, $\beta$ is the strength of the non-linearity, $\gamma$ is the
damping rate, $A$ is the forcing strength, $\omega$ is the external
forcing frequency, and $\omega_0$ is the natural frequency for small
amplitude oscillations in the undamped and unforced case
($\gamma=A=0$). By suitably normalising the time $t$ and position $x$,
and rescalling $\gamma, \omega$ and $A$, we can set $\omega_0=\beta=1$
without loss of generality and work with just three input
parameters. We shall use the term ``natural frequency'' to refer to
the frequency of the undamped and unforced oscillator. This frequency
depends on the amplitude of the oscillation $x_{max}$.  By mechanical
energy conservation, the natural frequency can be calculated from the
period of oscillation; $T(x_{max})=2\int_0^{x_{max}}dx/\dot{x}$, and
$\dot{x}^2/2=U(x_{max})-U(x)$, where the potential $U(x)$ is obtained
from $-U'(x)=-x-x^3$~\cite{Feynman:1963}.  In spite of being one of
the simplest nonlinear oscillators, the Duffing
oscillator~\eqref{eqn:duffing} exhibits surprisingly rich dynamical
behaviour owing to its amplitude-dependent natural frequency. More
precisely, as the forcing parameters $A$ and $\omega$ are varied, the
response can show symmetry-breaking pitchfork bifurcations,
period-doubling bifurcations, multistability (multiple attractors for
the same parameter settings) and even chaotic oscillations;
see~\cite{Holmes:1979,Parlitz:1985,Paar:1998,Kovacic:2011} for an
overview.

Whereas recent research has focused on nonlinear oscillators with
irregular external signals, there are aspects of
classical periodic forcing that have not been fully explored.  In this
paper, we relate the physical concept of subharmonic resonances and
mathematical concept of bifurcations to uncover a novel resonance
type. We obtain our main results for an intermediate damping strength
$\gamma=0.01$, and use bifurcation continuation
techniques~\cite{Doedel:2009} to calculate the structure of resonance
tongues in the plane of the the forcing frequency $\omega$ and
amplitude $A$. Furthermore, we describe the difference between
resonance and synchronisation~\cite{Pikovsky:2003}, and show that the
novel resonance is unusual in the sense that is has some of the
characteristics of synchronised oscillations.  Our results explain the
``patchy' tongues'' observed in recent numerical
simulations~\cite{Paar:1998} and link them to previous results on
instabilities and chaos due to ``bifurcation superstructure'' in
Eq.~\eqref{eqn:duffing}~\cite{Parlitz:1985}.

In order to define synchronisation, we consider a generalisation of
the damped Duffing oscillator~(\ref{eqn:duffing}) to the
self-sustained Duffing-Van der Pol oscillator:
\begin{equation}
\frac{d^2x}{dt^2} + \gamma(x)\frac{dx}{dt}  + x + x^3 = A \cos(\omega t),
\label{eqn:duffingvdp}
\end{equation}
where the damping rate $\gamma$ can depend on $x$ and change sign with
$x$, thus also act as a source of energy into the oscillator when its
amplitude becomes small.

\subsection{Bifurcation superstructure: Odd and even resonances}
Due to the nonlinear term $x^3$, which gives rise to amplitude
dependent natural frequency, Eq~\eqref{eqn:duffing} exhibits many more
resonances in addition to the main harmonic resonance near $\omega=1$.
For sufficiently small $A$ and $\gamma$, these resonances occur near
$\omega= 1/k$. They were first examined in detail by Parlitz and
Lauterborn~\cite{Parlitz:1985} (PL85), who refer to a resonance as odd
(even) when $k$ is odd (even).  When $A$ is increased, the resonances
shift away from $\omega= 1/k$ due to the forcing-induced change in the
oscillation amplitude and the resulting shift in the corresponding
natural frequency. (PL85) focus on moderate damping rate
($\gamma=0.2$) and high forcing strength ($0<A<50$), and perform
numerical bifurcation analysis which reveals self-similar set of
bifurcations, referred to as ``bifurcation superstructure'', with
regions of chaos in the $(\omega, A)$ parameter plane.  They associate
the ``bifurcation superstructure'' with the alternating odd and even
resonances.

\subsection{Patchy tongues}
Paar and Pavin~\cite{Paar:1998} (PP98) discuss coexisting attractors
in Eq.~\eqref{eqn:duffing} with weak damping ($\gamma=0.001$) and
moderate forcing strength ($0<A<5$) using numerical simulations.  This
parameter range is devoid of bifurcations giving rise to chaotic
oscillations. Instead, (PP98) demonstrate a high degree of
multistability between periodic solutions and show an intriguing
pattern of patchy tongues in the $(\omega, A)$ parameter plane; see
Fig. \ref{fig:duff_tongues} which was obtained in the same way
as~\cite[Fig.1]{Paar:1998} but for $\gamma=0.01$.  The colours in
Fig. \ref{fig:duff_tongues} refer to different integer ratios of the
period of periodic responses and the period $2\pi/\omega$ of the
forcing. This pattern cannot be explained in terms of even and odd
resonances identified in (PL85), and is referred to by (PP98) as
``intermingled Arnold tongues''.

In the following analysis we show that the backbone of the intriguing
pattern of intermingled tongues from (PP98) is formed by resonance
tongues associated with the novel resonance type and not by Arnold
tongues. Furthermore, we continue the largest resonance tongue of the
novel type from $\gamma=0.01$ to $\gamma=0.2$, and link it directly to
an element of the ``bifurcation superstructure'' from (PL85).

\section{Resonance vs. Synchronisation}
\label{Resonance vs Synchronisation}
In order to illustrate the main differences between the phenomenon of
resonance and the phenomenon of synchronisation, we compare in
Fig.~\ref{fig:comparison} the structure of the solutions to the
underdamped and forced Duffing oscillator~\eqref{eqn:duffing} with
$\gamma=0.01$ in panels (a) and (c), and the self-sustained and forced
Duffing-Van der Pol oscillator~\eqref{eqn:duffingvdp} with $\gamma(x)=
(1-x^2)\gamma$ in panels (b) and (d).  Following are definitions and
the listing of the characteristic properties for the two phenomena.

\subsection{Resonance}
A \emph{resonance} phenomenon is defined for a conservative oscillator
($\gamma(x)=0$) or a \emph{dissipative oscillator} [$\gamma(x)>0$ for
some $x$] as a noticeable increase in the amplitude of oscillations
near certain forcing frequencies.  This is best illustrated by fixing
$A$ and plotting the \emph{resonance curve}: the amplitude of
periodic solutions $x_{max}$ as a function of the forcing frequency
$\omega$.  Classical resonance in damped oscillators has the following
properties:
\begin{itemize}
\item (r1) In a linear oscillator (i.e.  without the $x^3$ term in
  Eq.~\eqref{eqn:duffing}), there is only one resonance near the
  natural frequency.  In a nonlinear oscillator, depending on the type
  of the nonlinearity, there can be additional \emph{subharmonic
    resonances} when the forcing frequency approaches a fraction or a
  multiple of the natural frequency~\cite{Lauterborn:1976}.

\item (r2) Typically, the frequency of the resonant response is the
  same as the forcing frequency.  The frequency of a non-resonant
  response is also the same as the forcing frequency. (Note: this
  situation may change for large $A$~\cite{Parlitz:1985}.)

\item (r3) As defined above, resonance is, in general, a quantitative
  phenomenon. However, in damped oscillators with amplitude-dependent
  frequency, resonance may involve bistability and qualitative changes
  in the solutions namely saddle-node bifurcations of periodic solutions.
\end{itemize}

For example, Fig.~\ref{fig:comparison}(a) shows the resonance curve
for the Duffing oscillator~\eqref{eqn:duffing} with $A=0.05$ near the
main resonance. The resonance curve `leans over' so that there is a
range of frequencies with two stable periodic solutions, one of which
has a much larger amplitude and corresponds to a resonant response;
the two solutions also have different phases. Which of these two
solutions the system settles to depends on initial conditions.  This
bistable range is bounded by two saddle-node bifurcations $SN$ of
periodic solutions~\cite{Kuznetsov:1998}, marked with red diamonds in
Fig.~\ref{fig:comparison}(a). By varying the amplitude and the
frequency of the forcing, we have:

\begin{itemize}
\item (r4) In the $(\omega, A)$ parameter plane, the corresponding
  codimension-one saddle-node bifurcation curves form \emph{resonance
    tongues}.  A resonance tongue has a tip for $A$ small but nonzero
  and $\omega$ near the natural frequency or near a
  fraction/multiple of the natural frequency. Moreover, the tip
  corresponds to a codimension-two cusp bifurcation
  $C$~\cite{Kuznetsov:1998}, where the two branches of the saddle-node
  bifurcation curve meet in a tangency
  [Fig.~\ref{fig:comparison}(c)]. Subharmonic resonance tongues appear
  for larger $A$ than the harmonic resonance
  tongue~\cite{Parlitz:1985}.
\end{itemize}

Note that for large enough $A$, resonance tongues may shift in
frequency, have additional cusp points and may interact with other
bifurcations via special codimension-two bifurcation points. On the
other hand, for values of $A$ below the cusp point $C$, the resonance
curve does not show any bistability and closely resembles the
resonance curve of a damped harmonic oscillator.

\subsection{Synchronisation}
The phenomenon of \emph{m:n synchronisation}, where
$m,n\in\mathbb{N}$, is defined for a dissipative but
\emph{self-sustained oscillator} [e.g. $\gamma(x)=(1-x^2)\gamma$ in
Eq.~\eqref{eqn:duffingvdp}] as a stable and fixed-in-time relationship
between the phases of the forcing $\phi(t)=\omega t$ and of the
oscillator $\varphi(t)$:
$$
0<\left| m\phi(t) - n\varphi(t)\right| \le 2\pi\;\Rightarrow\;  \frac{\omega}{\Omega}=\frac{n}{m},
$$
where $\Omega$ is the frequency of the (forced) oscillator.
Classical synchronisation to periodic forcing has the following
properties:
\begin{itemize}
\item (s1) A nonlinear oscillator can synchronise to periodic external
  forcing at various forcing frequencies provided that $m\omega/n$ is
  sufficiently close to the natural frequency.

\item (s2) The frequency of the (periodic) synchronised response is
  $\Omega = m \omega/n$, which by (s1) is close to the natural frequency. However, an
  unsynchronised response is quasiperiodic, meaning that it has
  components at the natural frequency, the forcing frequency, and
  linear combinations of both frequencies.
\end{itemize}

Fig.~\ref{fig:comparison}(b) shows a range of forcing frequencies
where the oscillator~\eqref{eqn:duffingvdp} synchronises in the ratio
1:1 to periodic external forcing.  The solid curve is a branch of
stable periodic solutions corresponding to the synchronised response,
while the dashed curve is a branch of unstable periodic solutions. The
two branches meet and disappear in a saddle-node bifurcation of
periodic solutions at both ends of the synchronisation range [red
diamonds in Fig.~\ref{fig:comparison}(c)]. Outside the synchronisation
range, the incommensurate ratio of the oscillator frequency and the
forcing frequency leads to stable quasiperiodic solutions, which
appear as disconnected dots in Fig.~\ref{fig:comparison}(b). Although
resonance and synchronisation involve the same bifurcation type, there
are important differences between these two nonlinear phenomena:

\begin{itemize}
\item (s3) As defined above, synchronisation is a qualitative
  phenomenon because synchronisation-desynchronisation transitions
  correspond to saddle-node bifurcations $SN$ of periodic solutions
  [Fig.~\ref{fig:comparison}(b)].  However, in contrast to the
  resonance phenomenon, there is normally just one stable solution at
  each value of $\omega$, meaning that there is no bistability of
  periodic solutions for small $A$.  (Note: This may change for large
  $A$ where one expects complicated dynamics due to break-up of
  invariant tori~\cite{Ostlund:1983,Broer:1998} and potentially more
  than one stable solution.)

\item (s4) In the $(\omega, A)$ parameter plane, the corresponding
  saddle-node bifurcation curves form \emph{m:n synchronisation
    tongues} or \emph{Arnold tongues}~\cite{Boyland:1986,Norris:1993};
  see Fig.~\ref{fig:comparison}(d) for an example of a 1:1
  synchronisation tongue in the Duffing-Van der Pol
  oscillator~\eqref{eqn:duffingvdp}.  A \emph{m:n} synchronisation
  tongue has a tip at $A=0$ and $\omega=n\omega_0/m$, where $\omega_0$
  is the frequency of the unforced oscillations. In contrast to the
  resonance tongue, this tip does not correspond to a cusp
  bifurcation.
  \\
\end{itemize}

The synchronisation phenomenon is captured in the $(\omega, A)$
parameter plane by an infinite but countable number of Arnold tongues.
Unlike resonance tongues, Arnold tongues originate at $A=0$ and
$\omega=n\omega_0/m$, where $m,n\in\mathbb{N}$. As $A$ is increased,
the tongues widen and may shift in frequency but they do not overlap
until some critical value of $A_c>0$. Above $A_c$, overlapping tongues
indicate break-up of invariant tori through various mechanisms
including homoclinic tangencies between stable and unstable manifolds
of saddle-type periodic solutions like the one indicated with a dashed
curve in Fig.~\ref{fig:comparison}(b). What is more, for $A$ large
enough, Arnold tongues may shift in frequency and develop cusp points,
and will typically interact with other bifurcations at special
codimension-two bifurcation points~\cite{Ostlund:1983,Broer:1998}. The
model example of Arnold tongues is the circle
map~\cite{Boyland:1986,Ott:2002}.

The remainder of this paper focuses on the detailed structure of resonances
in the Duffing oscillator.

\section{Complicated resonance curve: Three resonance types}
\label{Three resonance types}

To explain the high degree of multistability and the intermingled
structure of patchy tongues observed in numerical simulations, we
investigate the resonance structure in Eq.\eqref{eqn:duffing}. Besides
the harmonic resonance depicted in Fig.\ref{fig:comparison}(a), there
are subharmonic resonances in accordance with property (r1).  The
resonance curve in Fig.~\ref{fig:amplitude_resp} for $A=3$ shows that
subharmonic resonances: occur for $\omega<1$, become distinct only for
larger $A$ than the one used in Fig.\ref{fig:comparison} (a) and, as
$A$ is increased, the first to appear are odd subharmonic
resonances. Using the notation adopted from (PL85) we denote by $R_k$
the (subharmonic) resonances occurring when $\omega= (\mbox{natural
  frequency})/k$.

\subsection{Conventional Resonances}

The periodic solutions represented in the resonance curve are found
using the numerical continuation techniques AUTO~\cite{Doedel:2009}, which
also allows detection of unstable periodic solutions. In
Fig.~\ref{fig:resonance_1_kind} we zoom in on the subharmonic
resonances from Fig.~\ref{fig:amplitude_resp}.  Here it is seen that
parts of the resonance curve are `punctuated' by intervals of unstable
periodic solutions.  There are two types of punctuation. Firstly, in
the odd resonances ($k>1$ and odd), the stable and unstable solutions
merge in saddle-node bifurcations $SN$ (diamonds in
Fig.~\ref{fig:resonance_1_kind}). Secondly, the unstable solutions
found between the odd resonances are bounded by pairs of pitchfork
bifurcations $P$ (squares in Fig.~\ref{fig:resonance_1_kind}). These
bifurcations give rise to pairs of stable symmetry-broken periodic
solutions, shown in red in Fig.~\ref{fig:resonance_2_kind}, which
appear as mirror imaged orbits $x\leftrightarrow -x$ in the
$(x,\dot{x})$ plane [insets in
Fig.~\ref{fig:resonance_2_kind}]~\cite{Holmes:1979,Parlitz:1985}.
These solutions correspond to even subharmonic resonances ($k>1$ and
even) and become distinct for $A$ larger than the odd subharmonic
resonances.  What is more, they themselves undergo saddle-node
bifurcations giving rise to regions of multistability of
symmetry-broken periodic solutions. 

Both bifurcation types, that is $SN$ and $P$, indicate qualitative
changes in the solutions associated with (subharmonic) resonances in
accordance with property (r3). Interestingly, the odd and even
resonances alone cannot explain the patchy tongue structure from
Fig.~\ref{fig:duff_tongues}. This means that there must be additional
periodic solutions which are not connected to the primary branch
[black curve in Fig.~\ref{fig:resonance_2_kind}].

\subsection{Isolated Resonances}

Figure~\ref{fig:bifurcation_k_0_01_A_3} shows a one-dimensional
bifurcation diagram obtained by direct time integration starting with
the same initial condition $(x(0),\dot{x}(0))=(0,0)$ for each value of
$\omega$.  A comparison with Fig.~\ref{fig:resonance_2_kind}, obtained
using continuation techniques, reveals stable periodic solutions that
belong to stable branches of the resonance curve from
Fig.~\ref{fig:resonance_2_kind}, as well as additional stable periodic
solutions that are not present in
Fig.~\ref{fig:resonance_2_kind}. These additional solutions are
symmetric period-3 oscillations shown in
Fig.\ref{fig:trajectories_isolas}. What is more, numerical  continuation
initiated from numerical data for these additional solutions reveals
that they form ``isolas'' disconnected from the primary branch of
periodic solutions [green branches are disconnected from the black and
red branches in Fig.\ref{fig:resonance_3_kind}].

We have now from Fig.~\ref{fig:resonance_1_kind} through
Fig.~\ref{fig:resonance_2_kind} to Fig.~\ref{fig:resonance_3_kind}
identified three different resonance types:\\
(i) \emph{Odd resonances} which occur on the branch of primary
(symmetric) periodic
solutions that oscillate at the forcing frequency (black),\\
(ii) \emph{Even resonances} due to symmetry-broken periodic
solutions which bifurcate off the primary branch and also oscillate at the forcing frequency (red),\\
(iii) \emph{Isolated resonances} due to isolas of periodic solutions
with greatly increased amplitude of oscillations (green). This novel
resonance is unusual in the sense that it has certain properties akin to
synchronised oscillations. Firstly, unlike $R_k$ resonances, the
frequency of oscillation is a fraction of the forcing
frequency. Secondly, the periodic solutions involved are bounded by
saddle-node bifurcations in a way which is reminiscent of
Fig.~\ref{fig:comparison}(b) rather than
Fig.~\ref{fig:comparison}(a). As a consequence, cusp points at low $A$
are not expected in the corresponding resonance tongues.

Altogether, the three resonance types contribute to a complicated
resonance curve in Fig.\ref{fig:resonance_3_kind}, consisting of
various subharmonic resonances and disconnected components,
which give rise to many regions of multistability.

\section{Resonance tongues}

The three resonance types which were identified in the complicated
resonance curve from Fig.\ref{fig:resonance_3_kind} can be defined in
terms of saddle-node (diamonds) and pitchfork (squares) bifurcation
points.  In the two-dimensional $(\omega,A)$ parameter plane, these
bifurcation points can be continued with AUTO to obtain the
corresponding bifurcation curves which are referred to as resonance
tongues [property (r4)]. Fig. \ref{fig:continuation_tongues} shows
the resonance tongues computed for all three types of
resonances. These include: (black) saddle-node bifurcations of the
primary branch of periodic solutions which correspond to odd
subharmonic resonances, (red) pitchfork bifurcations on the primary
branch of periodic solutions which correspond to even subharmonic
resonances, (blue) saddle-node bifurcations of the symmetry-broken
periodic solutions which give rise to multistability of even
resonances, and (green) saddle-node bifurcations bounding the isolas
which correspond to isolated resonances.

\subsection{Explaining the ``patchy tongue structure''}

\begin{figure}[h]
\includegraphics[width=\linewidth]{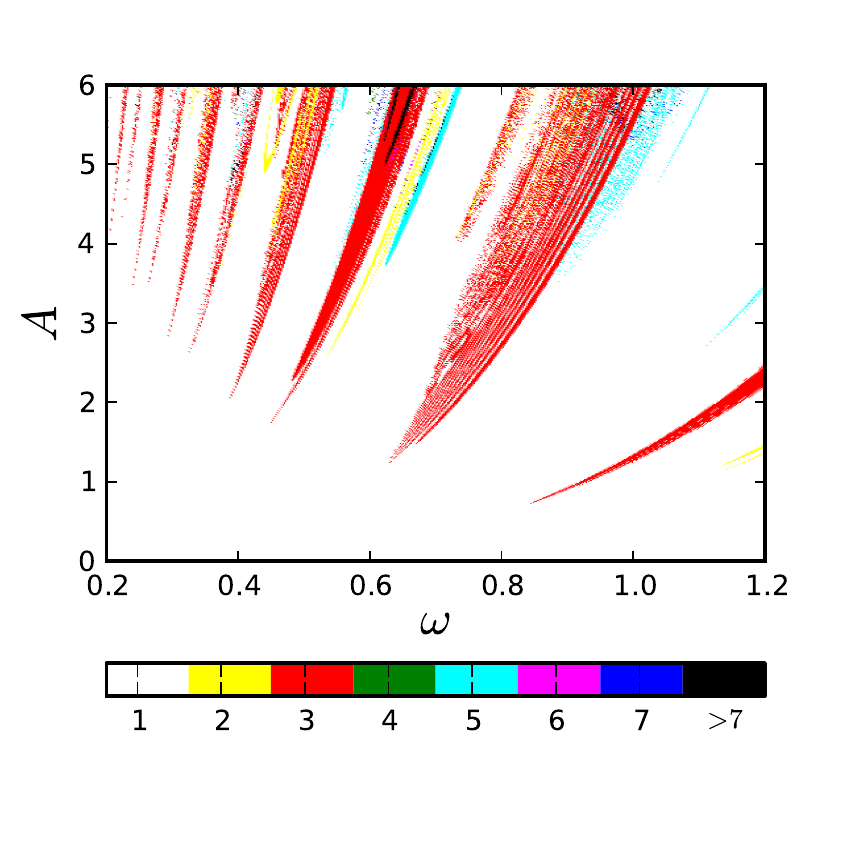}
\caption{Values of parameters $A$ and $\omega$ where we find
 attracting periodic orbits of the Duffing's system ~\eqref{eqn:duffing} of frequency (blank) $\omega$, (yellow)
 $\omega/2$, (red) $\omega/3$, (green) $\omega/4$, (cyan) $\omega/5$,
 (magenta) $\omega/6$, (blue) $\omega/7$, and (black) $\omega/n$ for $n=8,9,\ldots$. We used $\gamma=0.01$ and initial conditions $x(0) = \dot x(0) = 0$.}
\label{fig:duff_tongues}
\end{figure}

\begin{figure*}[h]
\includegraphics[width=6in, height=6in]{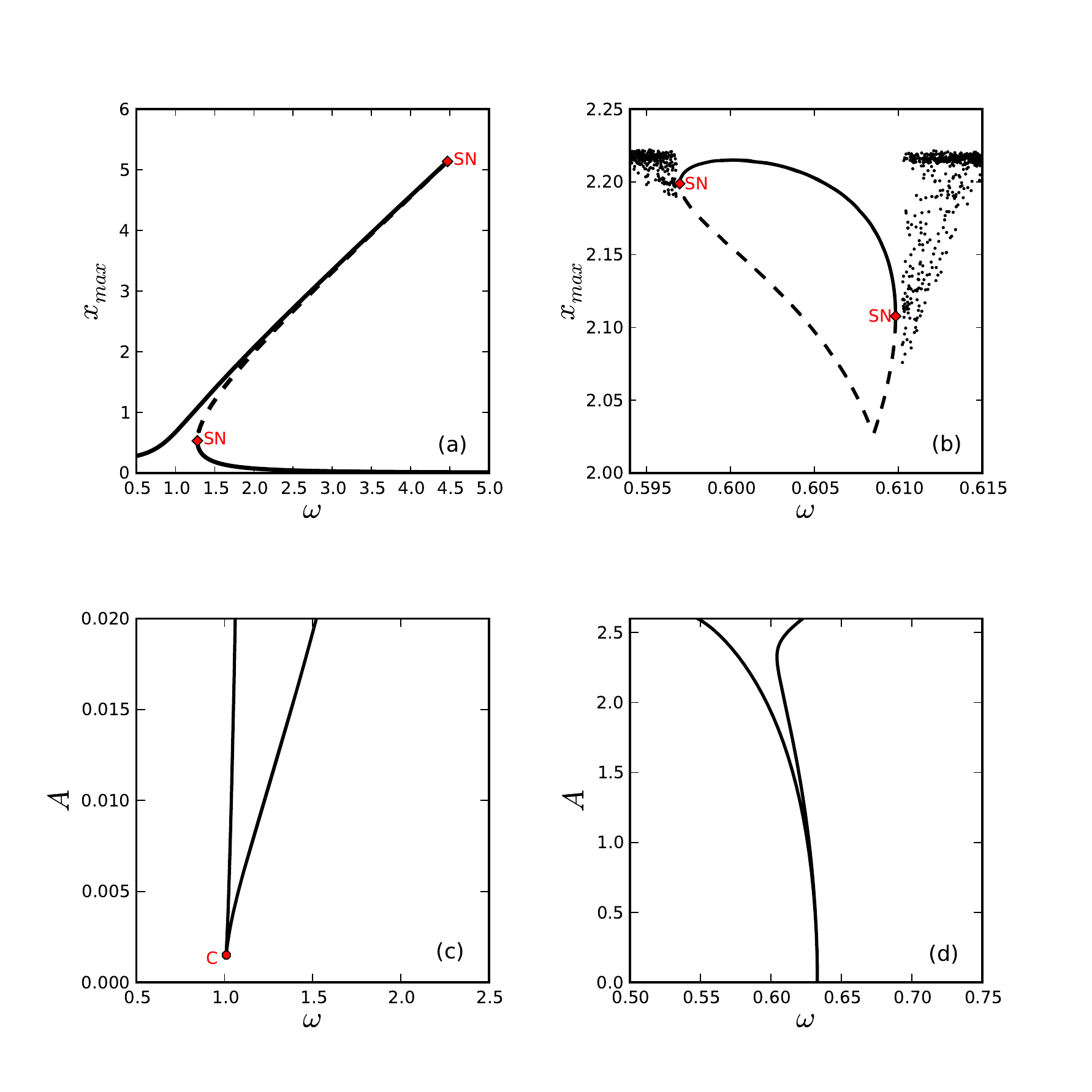}
\caption{(a) Amplitude resonance curve showing the primary resonance of the Duffing's system ~\eqref{eqn:duffing}. The black solid lines mark the stable solutions, the black dashed lines mark the unstable solutions connecting the two stable ones.  Saddle-node bifurcations (\emph{SN}) are marked by the red diamonds. Fixed parameters and initial condition: $A = 0.05$, $\beta= 1$, $\gamma = 0.01$; $x(0) = 0$, $\dot x(0)= 0$. (b) Values of the maximum of $x$ as a function of $\omega$ showing periodic and quasi-periodic oscillations of the Duffing-Van der Pol system \eqref{eqn:duffingvdp}. Periodic solutions are surrounded by quasi-periodic orbits (black dots). The black solid line marks the stable solutions, the black dashed line marks the unstable solutions. Fixed parameters and initial condition: $A = 2.0$, $\beta= 1$, $\gamma = 0.01$; $x(0) = \dot x(0) = 0 $. (c) Two parameter continuation ($\omega$, $A$) of the saddle-node bifurcations shown in panel (a). The cusp bifurcation \emph{C}, where the two saddle-node bifurcations merge, is marked by the red circle. Note that the cusp occur for $A>0$. (d) Two parameter continuation ($\omega$,$A$) of the saddle-node bifurcations shown in panel (b).}
\label{fig:comparison}
\end{figure*}

\begin{figure}[h]
\includegraphics[width=\linewidth]{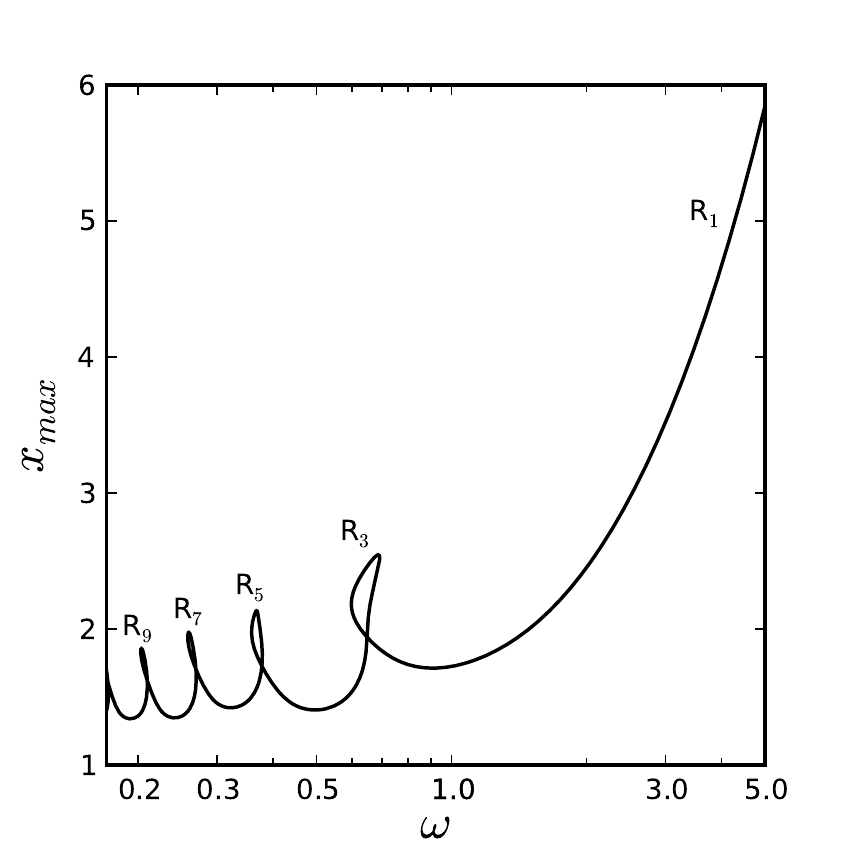}
\caption{ Amplitude resonance curve showing the maximum of the coordinate $x$ of the Duffing equation~\eqref{eqn:duffing} as a function of of the forcing frequency $\omega$. The curve is obtained with the numerical continuation method AUTO. Fixed parameters: $A = 3$, $\gamma = 0.01$. }
\label{fig:amplitude_resp}
\end{figure}

\begin{figure*}[h]
\includegraphics[width=18cm]{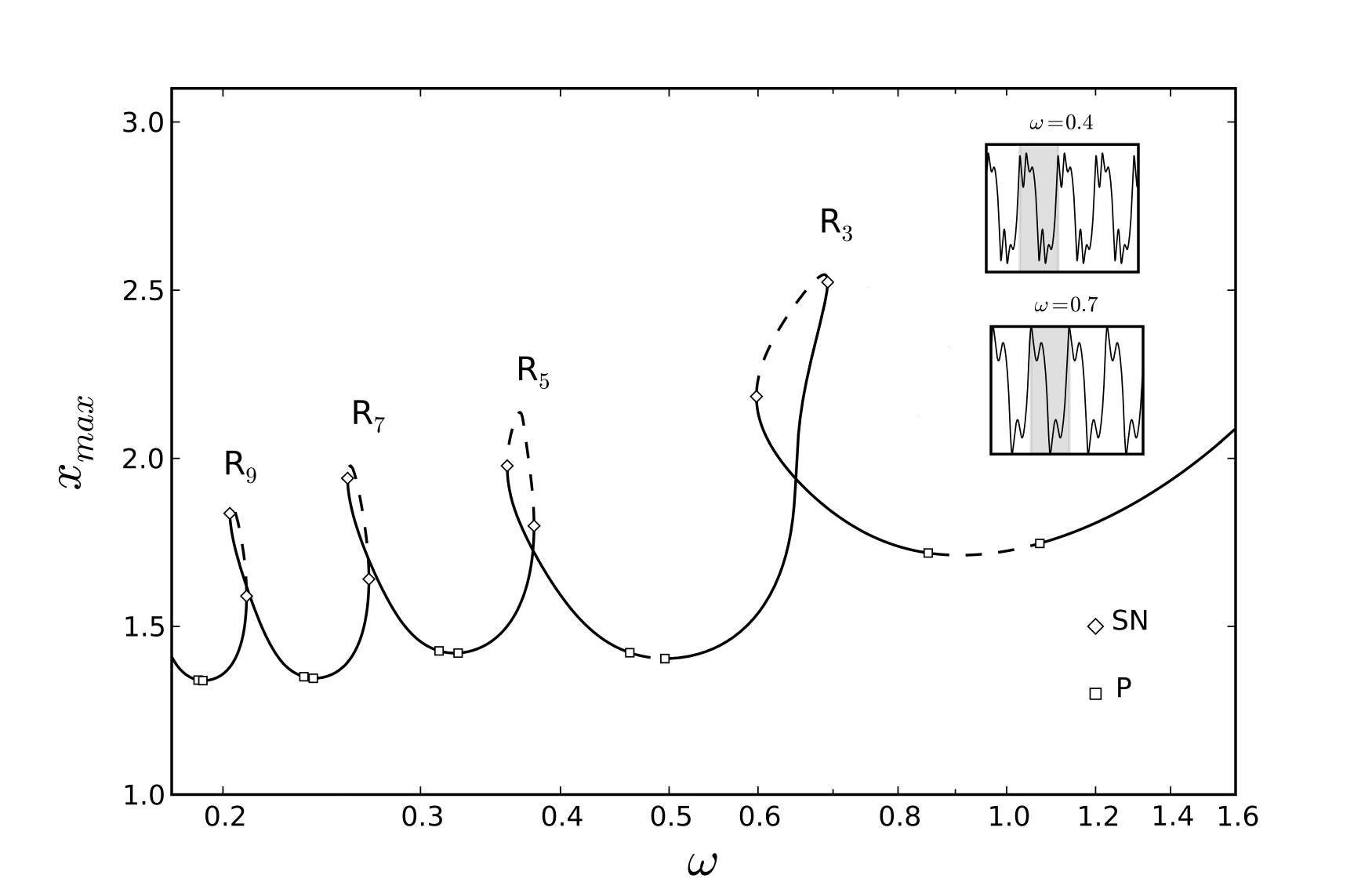}
\caption{Part of the resonance curve in Fig.\ref{fig:amplitude_resp} with branches of periodic solutions as a function of the forcing frequency $\omega$. On the primary branch of symmetric periodic solutions (black curve), saddle-node bifurcations (diamonds) and pitchfork bifurcations (squares) take place.  Stable branches are marked by solid lines, unstable branches by dashed lines. Fixed parameters: $A = 3 $,  $\gamma = 0.01$. The inserts show two solutions in the $(t,x)$ plane; for $\omega =0.4$ the solution is on the $R_5$ branch, where there are 5 local maxima in one period, marked by the gray band. For $\omega =0.7$ the solution is on the $R_3$ branch, where there are 3 local maxima in one period.}
\label{fig:resonance_1_kind}
\end{figure*}

\begin{figure*}[h]
\includegraphics[width=18cm]{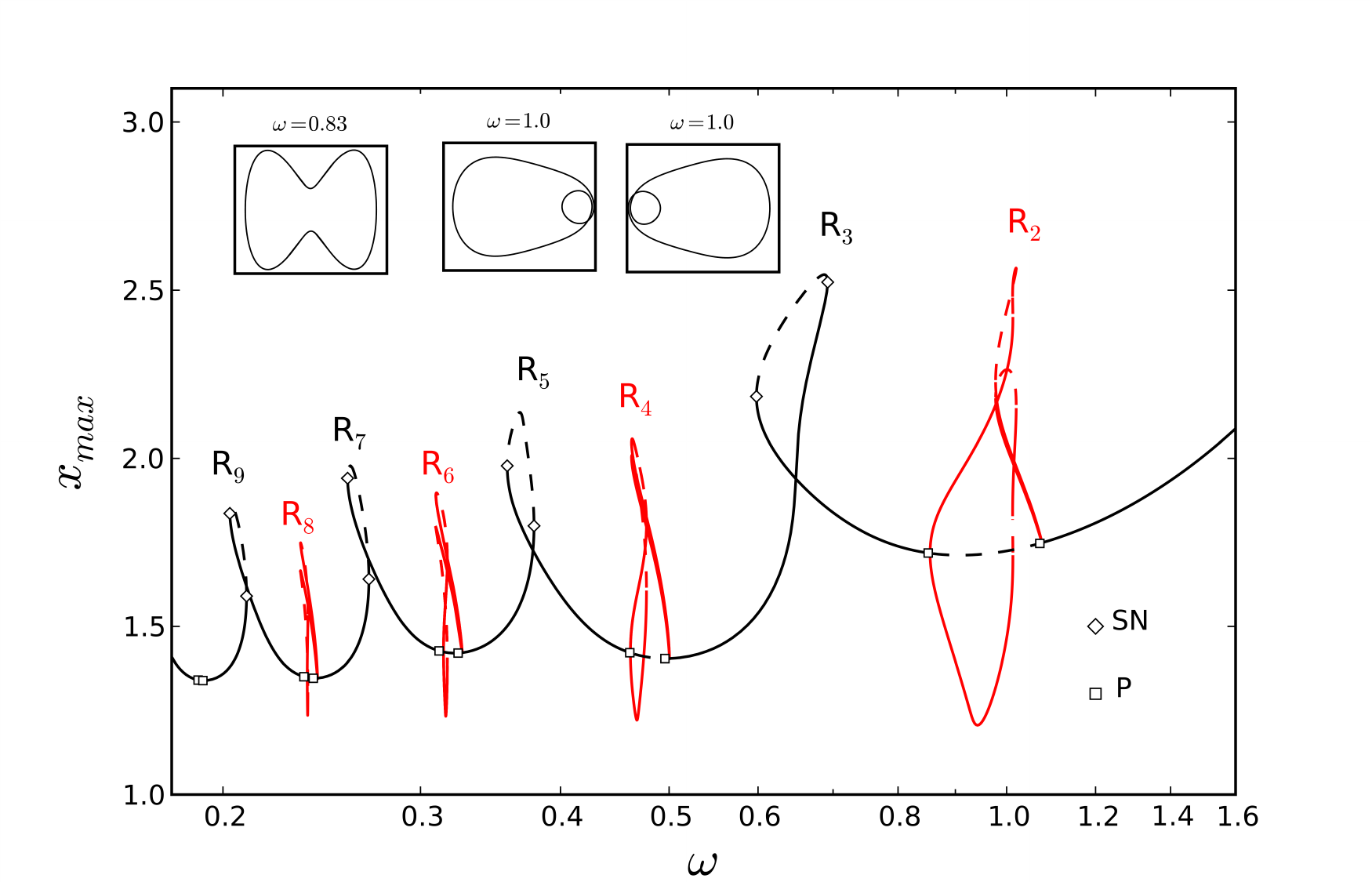}
\caption{ Same as Fig.\ref{fig:resonance_1_kind}, now including even resonances: on the primary branch of symmetric periodic solutions (black curve), saddle-node bifurcations (diamonds) and pitchfork bifurcations (squares) take place. At pitchfork bifurcations, pairs of stable symmetry-broken solutions appear (red curves). Stable branches are marked by solid lines, unstable branches by dashed lines.  Fixed parameters: $A = 3$,  $\gamma = 0.01$. Inserts show phase space trajectories $(x, \dot{x})$ just before a pitchfork bifurcation ($\omega=0.83$) and
a pair of solutions after the $x\rightarrow -x$ symmetry breaking ($\omega=1.0$). }
\label{fig:resonance_2_kind}
\end{figure*}

\begin{figure}[h]
\includegraphics[width=\linewidth]{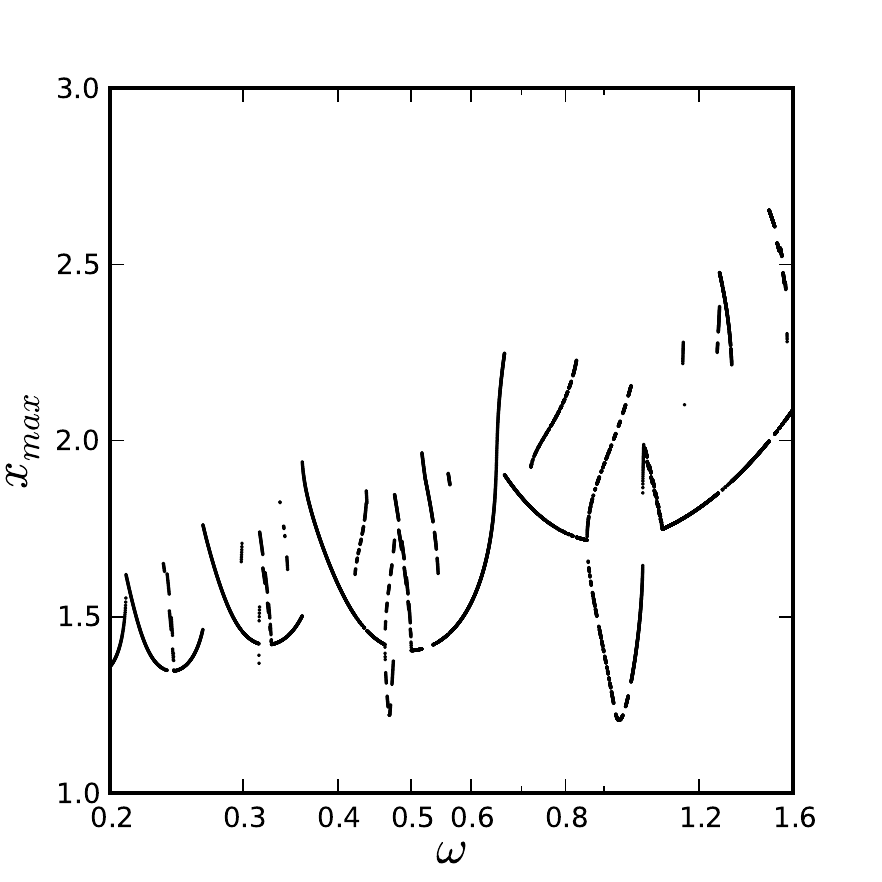}
\caption{1D numerical bifurcation diagram showing the maximum of $x$ plotted versus $\omega$ for $\gamma=0.01$, $A= 3$, $x(0) =  \dot x(0) = 0 $. }
\label{fig:bifurcation_k_0_01_A_3}
\end{figure}

\begin{figure*}[h]
\includegraphics[width=6in, height=6in]{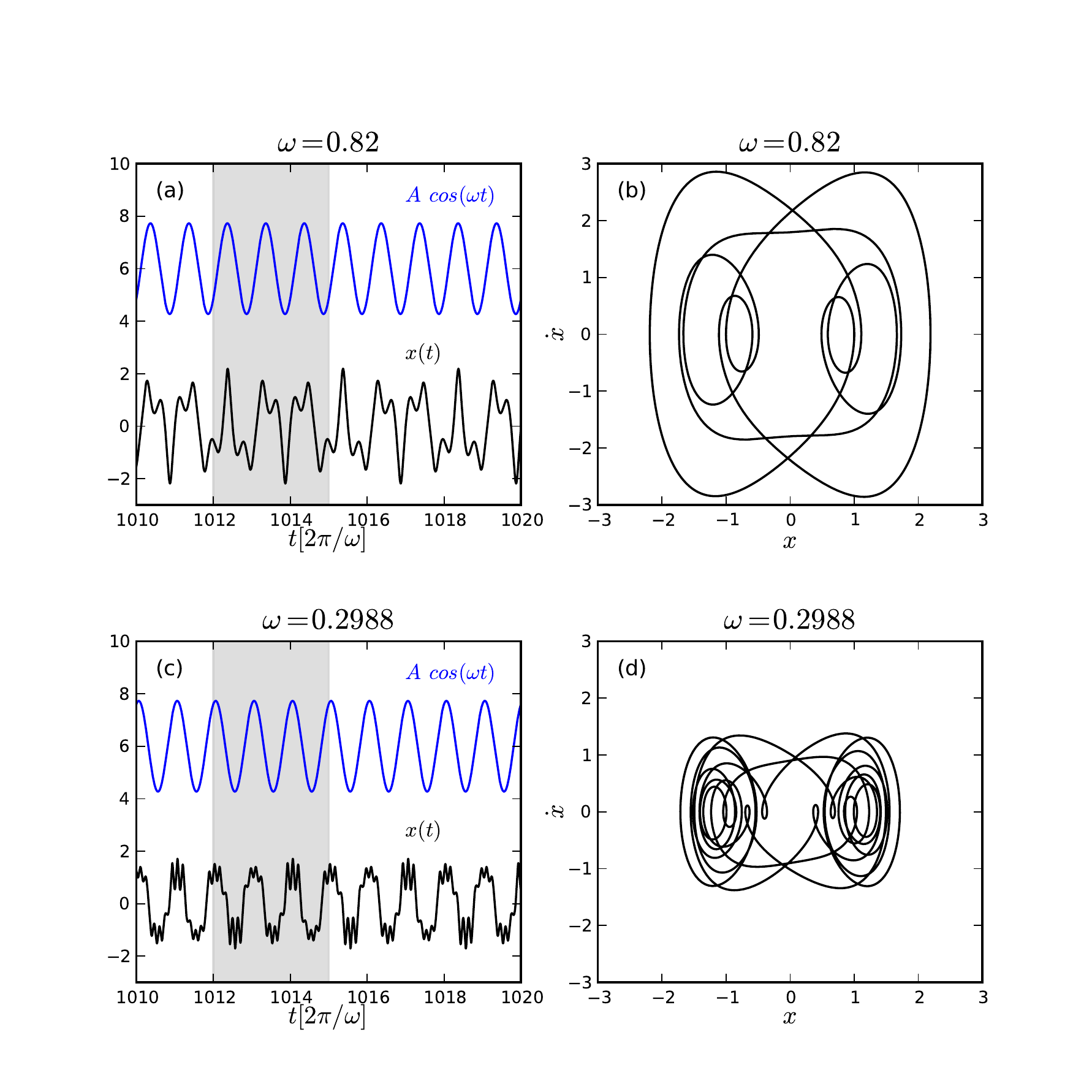}
\caption{Equation \eqref{eqn:duffing} has been integrated over 1020 periods of the forcing. The additional attracting orbits of the Duffing's system are shown for the last 10 periods of the forcing in the $(t,x)$ and $(x,\dot {x})$ planes. The solutions are period-3 orbits indicated by the gray bands. Fixed parameters and initial condition: $A = 3$, $\gamma= 0.01$, $x(0) =\dot x(0) = 0 $. (a) Trajectory in the $(t,x)$ plane , $\omega = 0.82$. (b) Phase space portrait in the two-dimensional $(x, \dot x)$ phase space, $\omega = 0.82$. (c) Trajectory in the $(t,x)$ plane, $\omega = 0.2988$. (d) Phase space portrait in the two-dimensional $(x, \dot x)$ phase space, $\omega = 0.2988$.}
\label{fig:trajectories_isolas}
\end{figure*}

\begin{figure*}[h]
\includegraphics[width=16cm]{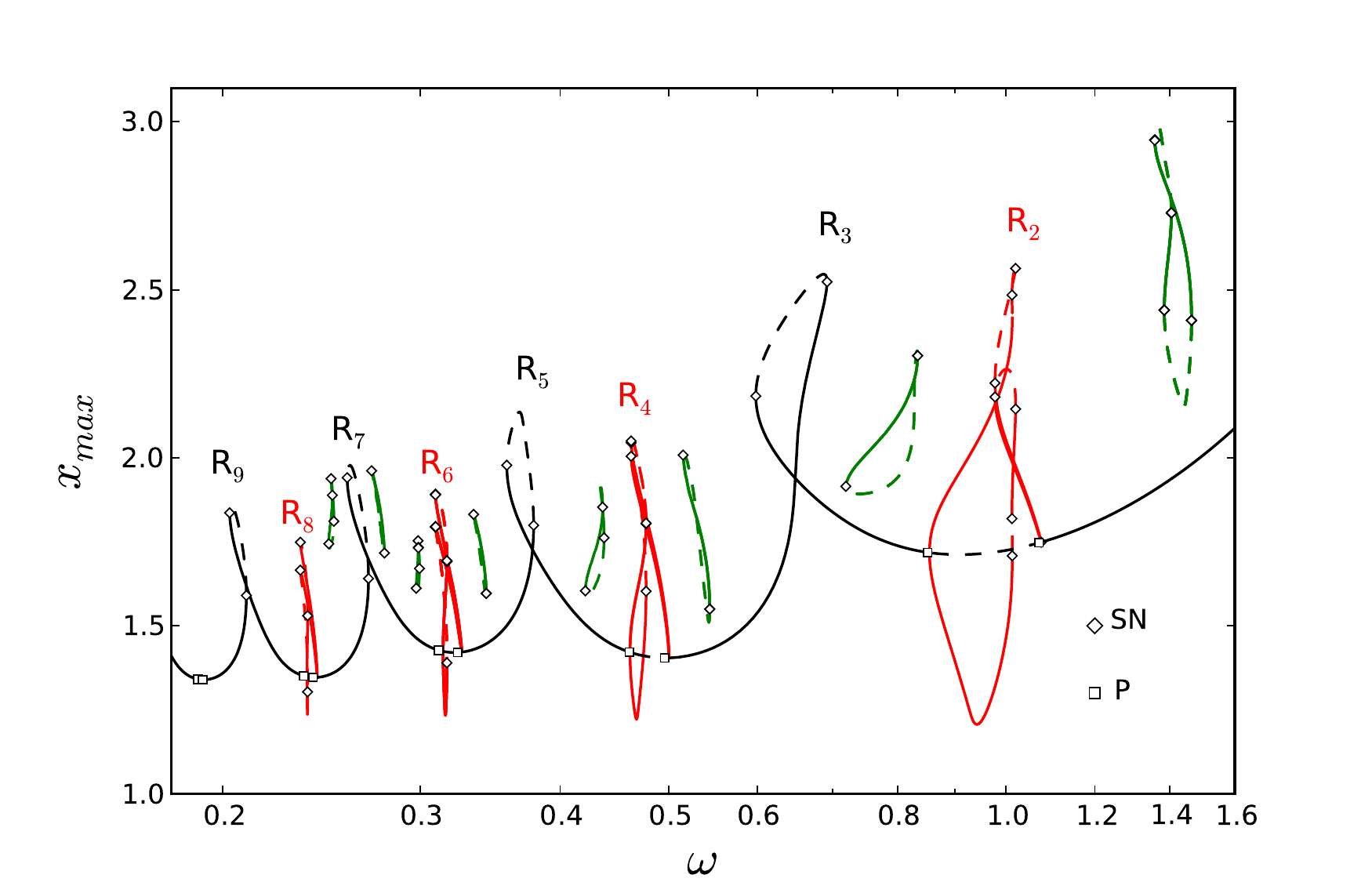}
\caption{Same 1D bifurcation diagram as Figures \ref{fig:resonance_1_kind} and \ref{fig:resonance_2_kind}, but now showing three different types of resonances: Resonances occurring on the symmetric branch of periodic solutions (black curve), resonances occurring on the asymmetric branch of periodic solutions (red curves), isolated resonances (green curves). Stable branches are marked by solid lines, unstable branches by dashed lines.  Fixed parameters: $A=3$, $\gamma = 0.01$.  }
\label{fig:resonance_3_kind}
\end{figure*}

\begin{figure*}[h]
\includegraphics[width=18cm]{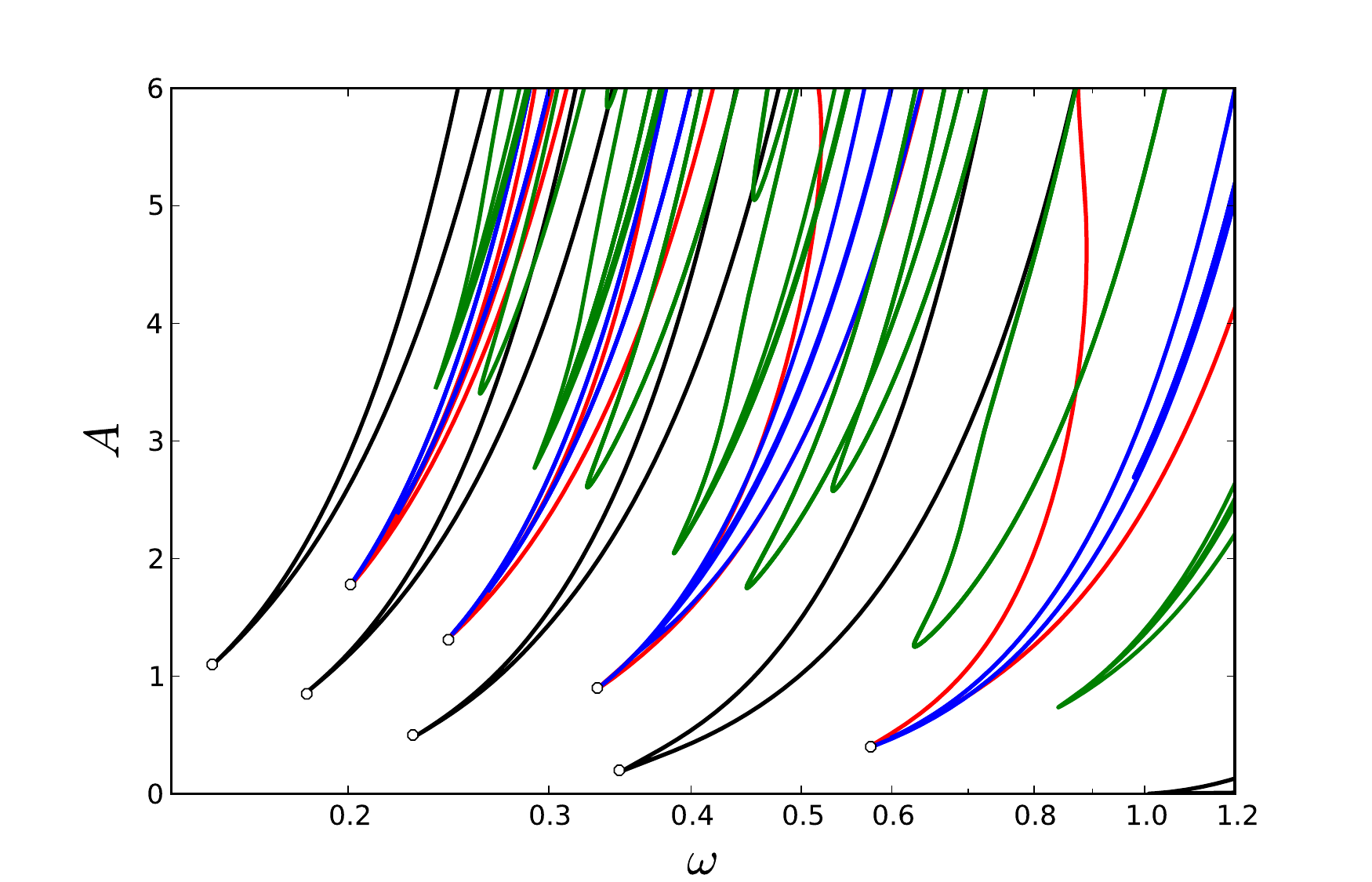}
\caption{ Resonance tongues computed
for all three types of resonances. These include: (black)
saddle-node bifurcations of the primary branch of periodic solutions which correspond to odd subharmonic
resonances, (red) pitchfork bifurcations on the primary
branch of periodic solutions which correspond to even
subharmonic resonances, (blue) saddle-node bifurcations
of the symmetry-broken periodic solutions which give rise
to multistability of even resonances, and (green) saddle-
node bifurcations bounding the isolas which correspond
to isolated resonances.  Using the notation adopted in (PL85), the resonance tongues computed for the resonance type (i) (\emph{odd resonances}) are denoted by $R_{1,1}$, $R_{3,1}$, $R_{5,1}$, $R_{7,1}$, $R_{9,1}$, and the resonance tongues computed for the resonance type (ii) (\emph{even resonances}) are denoted by $R_{2,1}$, $R_{4,1}$, $R_{6,1}$, $R_{8,1}$. The first subscript indicates the winding number (here defined as the number of maxima or minima of the periodic solution in one period) and the second subscript is the period in unit of the forcing period $2\pi/\omega$. The open circles indicate cusp bifurcations. Fixed parameter: $\gamma = 0.01$. 
}
\label{fig:continuation_tongues}
\end{figure*}

\begin{figure*}[h]
\includegraphics[width=18cm]{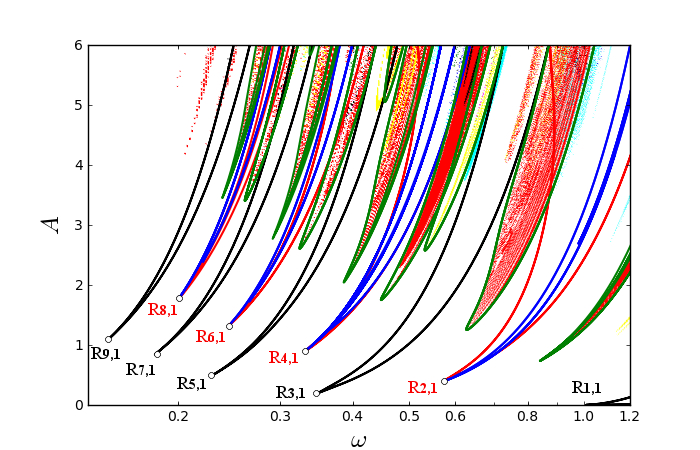}
\caption{Numerical estimates of parameters $A$ and $\omega$ leading to periodic attracting trajectories of period $nT_F$, $n = 1,...,7$ (same as Fig.\ref{fig:duff_tongues}) and the bifurcation boundaries of the resonance tongues obtained with the numerical continuation method AUTO are superimposed. The (red) regions with stable period-3 solutions and the (green) resonance tongues corresponding to the isolated resonances match perfectly. The open circles indicate cusp bifurcations. 
}
\label{fig:tongues_plus_continuation}
\end{figure*}

The structure of the resonance tongues bears strong resemblance to the
``intermingled Arnold tongues" reported by (PP98)
(Fig. \ref{fig:duff_tongues}). We can now superimpose the periodic
solutions oscillating at various fractions of the forcing frequency
from Fig. \ref{fig:duff_tongues} obtained by direct time integration
and the resonance tongues from Fig. \ref{fig:continuation_tongues}
obtained by continuation.  The resulting
Fig.\ref{fig:tongues_plus_continuation} reveals perfect match between
the (green) resonance tongues corresponding to the isolated resonances
and (red) regions wit stable periodic solutions oscillating at a
third of the forcing frequency. 
Thus, the intriguing patchy tongue structure found in (PP98) and shown
in Fig.~\ref{fig:duff_tongues} can be identified with the isolated
resonances, which appear to form the backbone of the structure.
Furthermore, the patchiness of the tongues is a result of
multistability: there are multiple stable periodic solutions of
different frequency for the same parameter settings and, as the
initial condition is fixed but the parameters $\omega$ and $A$ are
varied, the system can settle to a different periodic solution.

There are additional isolated resonances with stable periodic solutions
oscillating at $\omega/n$ for $n=2,4,5,\ldots$, whose resonance
tongues match the remaining patchy tongues from
Fig. \ref{fig:duff_tongues}.  For clarity, these additional isolated
resonance tongues are left out in Fig.\ref{fig:tongues_plus_continuation}.

\subsection{Links to ``bifurcation superstructure''}
\label{Comparison with previous works}

Our analysis of the intermediate damping rate $\gamma=0.01$ and
moderate forcing strength $0<A<6$ reveals various periodic solutions
and their bifurcations. However, we have not found any bifurcations,
such us torus bifurcations or period-doubling cascades, that would
eventually lead to irregular or chaotic oscillations. Rather, the
saddle-node and pitchfork bifurcations give rise to regions of
multistability of periodic solutions.

As the damping rate $\gamma$ is decreased, the patchy tongues become
more abundant. This was demonstrated by (PP98) for
$\gamma=0.001$. Hence, isolated resonances become even more prominent
at low damping, where they appear at even lower $A$ and give rise to a
greater degree of multistability of periodic solutions.

On the other hand, as the damping rate is increased, isolated
resonances seem to give way to the ``bifurcation superstructure''
involving period-doubling cascades, period-3 solutions and chaotic
attractors, which were reported for $\gamma=0.2$ and $10<A<50$ by
(PL85).

\begin{figure*}[h]
\includegraphics[width=6in, height=6in]{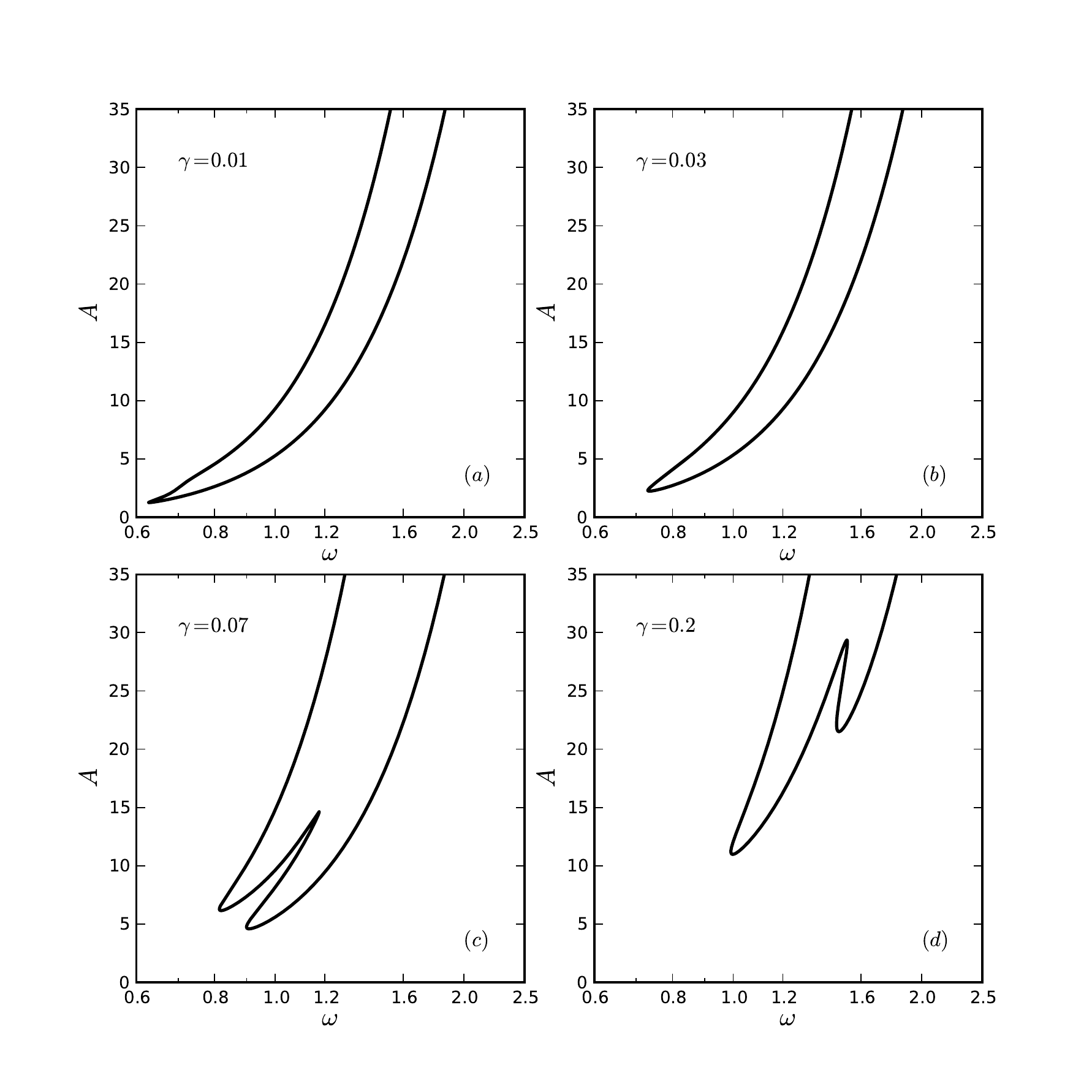}
\caption{ Bifurcation curve of the largest isolated resonance for different values of the damping parameter $\gamma$.}
\label{fig:largest_isola_different_k}
\end{figure*}

An interesting question emerges whether isolated resonances are purely
a low-damping phenomenon or whether the period-3 solutions associated
with isolated resonances (Fig.~\ref{fig:trajectories_isolas}) are
related to the ``bifurcation superstructure''
from~\cite{Parlitz:1985}.  To address this question we continued the
largest isolated resonance tongue from
Fig.\ref{fig:continuation_tongues} to higher values of $\gamma$. The
results of the continuation are shown in
Fig.\ref{fig:largest_isola_different_k}.  As $\gamma$ is increased
from $\gamma=0.01$ [Fig.\ref{fig:largest_isola_different_k}(a)], the
isolated resonance tongue moves to higher values of $A$ and develops
another tip [Fig.\ref{fig:largest_isola_different_k}(b-c)]. At
$\gamma=0.2$ the continuation of the isolated resonance tongue matches
exactly an element of the ``bifurcation superstructure'' denoted with
$R_{9,3}$ in~\cite[Fig.6]{Parlitz:1985}.  Thus, in addition to
explaining ``patchy tongues'', isolated resonances (i) provide a link
between ``patchy tongues'' the ``bifurcation superstructure'' and (ii) give
new insight into ``bifurcation superstructure''~\cite{Parlitz:1985}
which now appears to be organised by three rather than two different
resonance types.

\begin{figure}[h]
\includegraphics[width=\linewidth]{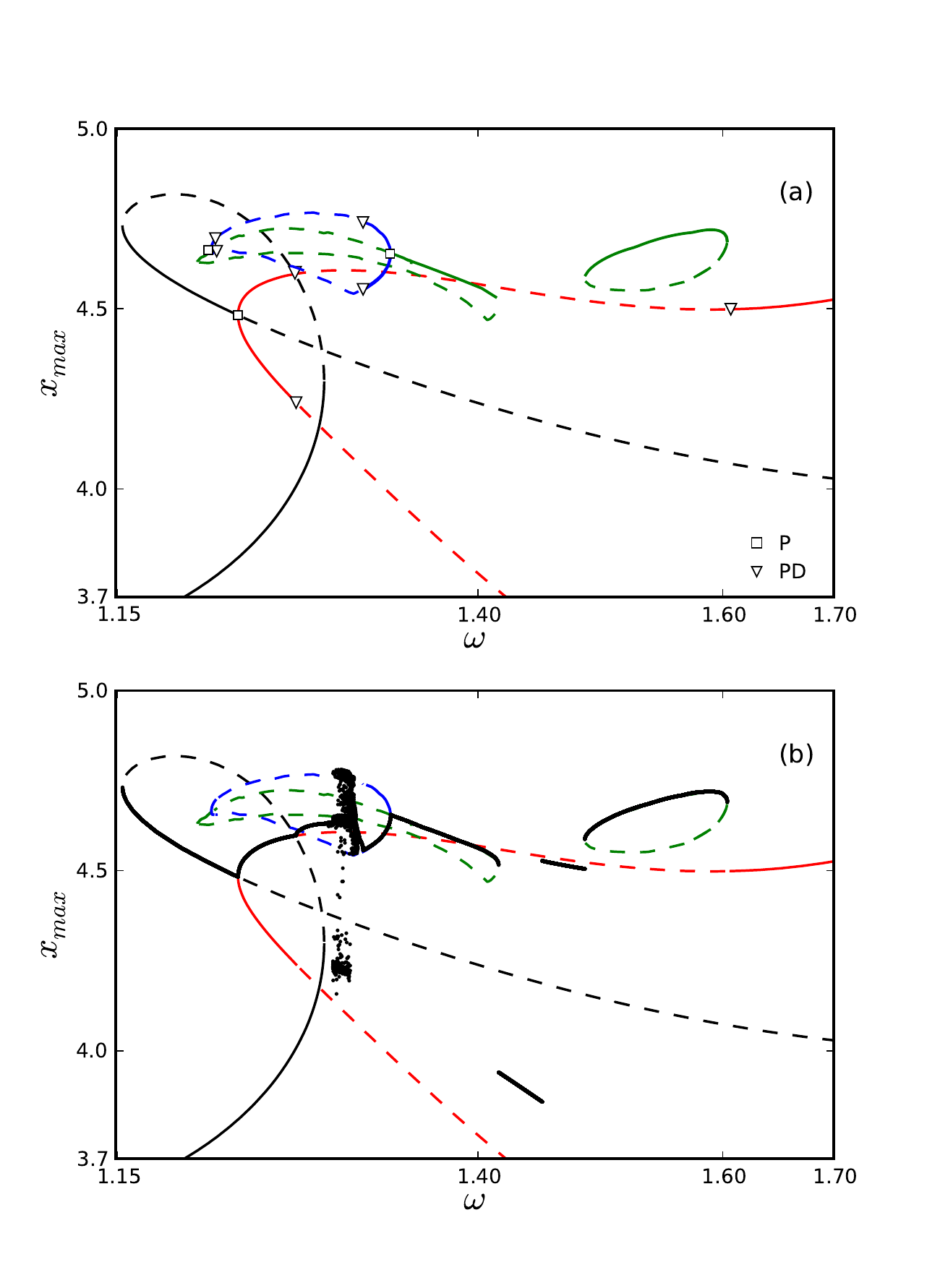}
\caption{(a) AUTO bifurcation diagram in the range $1.15 < \omega <  1.7$  showing the primary branch of periodic solutions (black), symmetry-broken solutions bifurcating off the primary branch (red), isolas of periodic solutions (green), and symmetry-broken solutions bifurcating from the left isola (blue). (b) A numerical bifurcation diagram is superimposed on the AUTO continuation. The numerical bifurcation diagram consists of four runs: starting at the stable right and left isolas of periodic solutions, the control parameter is increased and decreased.    }
\label{fig:bifurcation_A_25_k_0_2_duff}
\end{figure}

The key difference is that periodic solutions involved in isolated
resonances at low $\gamma$ may cease to be isolated as $\gamma$ is
increased.  The bifurcation diagram in
Fig.~\ref{fig:bifurcation_A_25_k_0_2_duff}(a), obtained using
continuation techniques~\cite{Doedel:2009} for $\gamma = 0.2$, $A=25$
and $1.15 \leq \omega \leq 1.7$, shows the primary branch of periodic
solutions (black), the symmetry-broken solutions bifurcating off the
primary branch (red) via pitchfork bifurcation (square), two
potentially isolated components of periodic solutions (green), and
symmetry-broken solutions (blue) bifurcating from the left potentially
isolated component via pitchfork bifurcation (square).  In
Fig.~\ref{fig:bifurcation_A_25_k_0_2_duff}(b), a numerical bifurcation
diagram obtained starting on the stable branch of each (green)
potentially isolated component is superimposed on the continuation
results from panel (a). The right-hand side component turns out to be
isolated. However, the larger left-hand side component is connected to
other solutions. Starting at the stable periodic solution of this
component and decreasing $\omega$, the solution looses stability via
symmetry-breaking pitchfork bifurcation, then the symmetry-broken
solution undergoes period-doubling cascade to chaos, followed by an
inverse period-doubling cascade to period-one solution (red) which, in
turn, connects to the primary branch of periodic solutions (black) via
pitchfork bifurcation.

\section{Conclusion}

We have investigated nonlinear resonances in terms of
bifurcations of periodic solutions using the example of a periodically
forced Duffing oscillator. We have identified novel isolated
resonances in addition to already known odd and even subharmonic resonances, and
demonstrated a complicated resonance `curve' with isolas (isolated
components) of periodic solutions and high degree of
multistability. Most importantly, the identified here isolated
resonances in conjunction with numerical continuation techniques
allowed us to (i) explain the intriguing structure of ``patchy
tongues'' observed in stability diagrams for weak damping and (ii)
link those ``patchy tongues'' to a seemingly unrelated phenomenon of
``bifurcation superstructure'' found for moderate damping.

Firstly, to avoid confusion between resonance and synchronisation, we
have defined each phenomenon in terms of bifurcations, gave a short
discussion of the key differences, and described how the two phenomena
exhibit themselves in the stability diagrams. In particular, we have
distinguished between resonance tongues and synchronisation tongues,
which are also known as Arnold tongues, in the parameter plane of the
forcing frequency vs. the forcing amplitude.  Secondly, we have shown
that resonance tongues associated with the isolated resonances form
the backbone of the intriguing pattern of ``patchy tongues" at weak
damping, which were reported by Paar and Pavin and mistaken for Arnold
tongues~\cite{Paar:1998}. Thirdly, we have demonstrated that, as the
damping rate increases, these resonances may cease to be isolated and
the corresponding resonance tongues evolve into particular elements of
the ``bifurcation superstructure" reported by Parlitz and Lauterborn
for moderate damping~\cite{Parlitz:1985}.

The new insight into nonlinear resonances in the simple Duffing
oscillator can be extended to more complicated systems. One example
are class-B lasers which exhibit damped oscillations (relaxation
oscillation onto the lasing solution) and self-sustained oscillations
(the lasing solution itself) at the same time. Stability diagrams for
lasers subject to external optical
injection~\cite[Figs. 9-11]{Krauskopf:2002} show variety of coexisting
tongues, which can be interpreted as a combination of a 1:1
synchronisation tongue and different types of nonlinear resonance tongues.

\bibliography{qhe}

\end{document}